\documentclass[a4paper,10pt]{amsart}
\usepackage{pst-all, pstricks, pst-node}
\usepackage{subfigure}
\usepackage{graphicx, amssymb, amsmath, amscd, hhline}
\usepackage[arrow,tips,matrix]{xy}

\theoremstyle{ams}
\newtheorem{theorem}{Theorem}[section]
\newtheorem{proposition}[theorem]{Proposition}
\newtheorem{lemma}[theorem]{Lemma}
\newtheorem{corollary}[theorem]{Corollary}

\theoremstyle{definition}
\newtheorem{question}[theorem]{Question}

\newtheorem{example}[theorem]{Example}

\newtheorem*{acknowledgement}{Acknowledgement}

\numberwithin{table}{section}
\numberwithin{figure}{section}
\numberwithin{equation}{section}

\newcommand{\R}{\mathbb{R}}

\newcommand{\SU}{\mathrm{SU}}
\newcommand{\su}{\mathfrak{su}}
\newcommand{\torus}{\mathfrak{t}}
\newcommand{\SO}{\mathrm{SO}}
\newcommand{\so}{\mathfrak{so}}

\begin{document}
\title[Self-indexing moment map]
{Hamiltonian circle action with self-indexing moment map}

\author[Y. Cho]{Yunhyung Cho}
\address{School of Mathematics, Korea Institute for Advanced Study,
87 Hoegiro, Dongdaemun-gu, Seoul, 130-722, Republic of Korea}
\email{yhcho@kias.re.kr}

\author[M. K. Kim]{Min Kyu Kim*}
\address{Department of Mathematics Education,
Gyeongin National University of Education, 45
Gyodae-Gil, Gyeyang-gu, Incheon, 407-753,
Republic of Korea} \email{mkkim@kias.re.kr}

\thanks{* the corresponding author}

\maketitle

\begin{abstract}
Let $(M,\omega)$ be a $2n$-dimensional smooth compact
symplectic manifold equipped with a Hamiltonian circle
action with only isolated fixed points and let $\mu : M
\rightarrow \R$ be a corresponding moment map. Let
$\Lambda_{2k}$ be the set of all fixed points of index
$2k$. In this paper, we will show that if $\mu$ is constant
on $\Lambda_{2k}$ for each $k$, then $(M,\omega)$ satisfies
the hard Lefschetz property. In particular, if $(M,\omega)$
admits a self-indexing moment map, i.e. $\mu(p) = 2k$ for
every $p \in \Lambda_{2k}$ and $k=0,1,\cdots,n,$ then
$(M,\omega)$ satisfies the hard Lefschetz property.
\end{abstract}

\section{Introduction} \label{section:
introduction}
Let $(M,\omega)$ be a $2n$-dimensional compact
symplectic manifold. We say that $\omega$
satisfies the \textit{hard Lefschetz property} if
\begin{displaymath}
\begin{array}{cccc}
\omega^{n-k} : & H^{k}(M) & \longrightarrow &
H^{2n-k}(M)\\ 
& \alpha & \mapsto & \alpha
\wedge \omega^{n-k}
\end{array}
\end{displaymath}
is an isomorphism for every $k=0, 1, \cdots, n$.
According to \cite{JHKLM}, the following question has been around for many years.

\begin{question} \label{conjecture} \cite{JHKLM}
Let $(M,\omega)$ be a closed symplectic manifold
with a Hamiltonian circle action. Assume that all
fixed points are isolated. Then does $(M,\omega)$
satisfy the hard Lefschetz property?
\end{question}

In this paper, we discuss about conditions under which a
Hamiltonian circle action with isolated fixed points
satisfies the hard Lefschetz property. A leading candidate
for these conditions might be `\textit{semifree}', i.e. the
action is free outside the fixed point set. We note that a
Hamiltonian circle action with isolated fixed points is
semifree if and only if the weights on the normal bundles
to the fixed points are all $\pm 1,$ hence `semifree' is a
condition on weights of $S^1$-representations on the fixed point set. Actually, a Hamiltonian semifree
circle action with isolated fixed points satisfies the hard
Lefschetz property as Reyer Sjamaar pointed out in
\cite{JHKLM}. But there is no other result related to
Question \ref{conjecture} as far as the authors know. A motivation of this paper is a resent result of the authors \cite{CK} which says that, we can determine if a
Hamiltonian torus action satisfies the hard Lefschetz property just by looking at its moments map image in some
cases. More precisely, the authors proved that a six-dimensional Hamiltonian GKM two-torus action with an
index increasing GKM graph satisfies the hard Lefschetz
property, see \cite{CK}. Hence one might guess that an image of a moment map, especially an image of \textit{1-skeleton} determines the hard Lefschetz property. In this point of
view, it is also conceivable to find such condition on the image of fixed point set in
the case of a Hamiltonian circle action. An elementary theory of Hamiltonian group action says that a moment map for a
Hamiltonian circle action with isolated fixed points is a
Morse function. A critical point set of a moment map is equal to
the fixed point set, see \cite{Au} for the details. In this paper, we prove
the followings.

\begin{theorem} \label{theorem: main}
Let $(M,\omega)$ be a $2n$-dimensional smooth compact
symplectic manifold equipped with a Hamiltonian circle
action with only isolated fixed points with a moment map
$\mu : M \rightarrow \R$. Let $\Lambda_{2k}$ be the set of
all fixed points of Morse index $2k$ with respect to $\mu.$
If $\mu$ is constant on $\Lambda_{2k}$ for each $k=0,1,
\cdots, n$, then $(M,\omega)$ satisfies the hard Lefschetz
property.
\end{theorem}

Note that the condition `\textit{constant $\mu$ on each $\Lambda_{2k}$}' is a
weaker version than the well known condition
`self-indexing' on Morse funtions. Morse function $f : M
\rightarrow \R$ is called \textit{self-indexing} if
$\lambda(p) = f(p)$ for every critical point $p \in M$,
where $\lambda(p)$ is a Morse index at $p,$ see \cite[p.
44, Definition 4.9]{Mi}. Hence if $(M,\omega)$ admits a
self-indexing moment map $\mu,$ then it satisfies the
condition in Theorem \ref{theorem: main} automatically. As
a corollary of Theorem \ref{theorem: main}, we have the
followings.

\begin{corollary}\label{self}
Let $(M,\omega)$ be a $2n$-dimensional smooth
compact symplectic manifold equipped with a
Hamiltonian circle action with only isolated
fixed points with a self-indexing moment map $\mu
: M \rightarrow \R$. Then $(M,\omega)$ satisfies
the hard Lefschetz property.
\end{corollary}

In Section \ref{section: example}, we give an
example of periodic Hamiltonian system which
admits a self-indexing moment map. This paper is
organized as follows. In Section \ref{section:
canonical class}, we will give a brief
introduction to equivariant cohomology for
Hamiltonian circle actions. In Section
\ref{section: proof}, we will prove Theorem
\ref{theorem: main}. In Section \ref{section:
example}, we will give two examples.

\section{Equivariant symplectic forms and
Canonical classes} \label{section: canonical
class}

In this section, we briefly review an elementary
equivariant cohomology theory for Hamiltonian circle
actions which will be used in the rest of the paper.
Throughout this section, we will assume that every
coefficient of any cohomology theory is $\R$. Let $S^1$ be
the unit circle group and let $M$ be an $S^1$-manifold.
Then the equivariant cohomology $H^*_{S^1}(M)$ is defined
by $H^*_{S^1}(M) := H^*(M \times_{S^1} ES^1)$ where $ES^1$
is a contractible space on which $S^1$ acts freely. Note
that $M \times_{S^1} ES^1$ has a natural $M$-bundle
structure over the classifying space $BS^1 := ES^1 / S^1$ so that
$H^*_{S^1}(M)$ admits a
$H^*(BS^1)$-module structure. For the fixed point set
$M^{S^1}$, the inclusion map $i : M^{S^1} \hookrightarrow
M$ induces a $H^*(BS^1)$-algebra homomorphism $$i^* :
H^*_{S^1}(M) \rightarrow H^*_{S^1}(M^{S^1}) \cong
\bigoplus_{F \subset M^{S^1}} H^*(F)\otimes H^*(BS^1)$$ and
we call $i^*$ \textit{the restriction map to the fixed
point set}. Note that for any fixed component $F \subset
M^{S^1}$, the inclusion map $i_F : F \hookrightarrow
M^{S^1}$ induces a natural projection $i_F^* :
H^*_{S^1}(M^{S^1}) \rightarrow H^*_{S^1}(F) \cong
H^*(F)\otimes H^*(BS^1).$ For every $\alpha \in
H^*_{S^1}(M)$, we will denote by $\alpha|_F$ an image
$i_F^* (i^*(\alpha)).$ Now, assume that $(M,\omega)$ be a smooth compact symplectic manifold with a Hamiltonian circle action with a moment map $\mu : M \rightarrow \R$. For each fixed component $F \subset M^{S^1}$, let $\nu_{F}$ be a normal
bundle of $F$ in $M$. Then the \textit{negative
normal bundle} $\nu^{-}_{F}$ of $F$ can be defined as
a sub-bundle of $\nu_F$ whose fiber over $p \in F$ is a subspace of $T_p M$ tangent to an unstable submanifold of $M$ at $F$ with respect to $\mu$ for every $p \in F$. We denote by $e^{-}_F \in H^*_{S^1}(F)$ the equivariant Euler class of $\nu^{-}_F$. In Hamiltonian cases, $H^*_{S^1}(M)$ has a remarkable property (called a \textit{Kirwan's injectivity theorem}) as
follows.

\begin{theorem} \cite{Ki} \label{theorem: formal}
Let $(M,\omega)$ be a closed symplectic manifold with a
Hamiltonian circle action. For an inclusion $f : M
\hookrightarrow M \times_{S^1} ES^1$ as a fiber, the
induced map $f^* : H^*_{S^1}(M) \rightarrow H^*(M)$ is
surjective. Equivalently, $H^*_{S^1}(M)$ is a free
$H^*(BS^1)$-module, i.e. $M$ is equivariantly formal.
Moreover, the kernel of $f^*$ is given by $u \cdot
H^*_{S^1}(M)$, where $u$ is a degree-two generator of $H^*(BS^1)$ and $\cdot$ is a scalar product as an $H^*(BS^1)$-module.
\end{theorem}

From now on, we assume that $(M,\omega)$ is a closed
symplectic manifold equipped with a Hamiltonian
$S^1$-action. By Theorem \ref{theorem: formal}, the
equivariant cohomology $H^*_{S^1}(M)$ is a free
$H^*(BS^1)$-module so that $H^*_{S^1}(M)$ is isomorphic to
$H^*(M) \otimes H^*(BS^1)$ as a $H^*(BS^1)$-module. Note
that $H^*(BS^1)$ is isomorphic to $\R[u]$ where $-u$ is the
first Chern class of the principal $S^1$-bundle $ES^1
\rightarrow BS^1$. D. McDuff and S. Tolman found a
remarkable family of equivariant cohomology classes as
follows.

\begin{theorem}\cite{McT}\label{McT}
Let $(M,\omega)$ be a closed symplectic manifold
equipped with a Hamiltonian circle action with a
moment map $\mu : M \rightarrow \R$. For each
connected fixed component $F \subset M^{S^1}$,
let $k_F$ be the index of $F$ with respect to
$\mu$.  Then given any cohomology class $Y \in
H^i(F)$, there exists a unique class
$\widetilde{Y} \in H_{S^1}^{i + k_F}(M)$ such
that
\begin{enumerate}
\item $\widetilde{Y}|_{F'} = 0$
for every $F' \in M^{S^1}$ with $\mu(F') <
\mu(F)$,
\item $\widetilde{Y}|_F = Y \cup e^{-}_F$, and
\item the degree of
$\widetilde{Y}|_{F'} \in H^*_{S^1}(F')$
is less than the index $k_{F'}$ of $F'$
for all fixed components $F' \neq F.$
\end{enumerate}
\end{theorem}

We call such a class $\widetilde{Y}$ \textit{the canonical
class} with respect to $Y$. In the case when all fixed
points are isolated, Theorem \ref{McT} implies the
following corollary.

\begin{corollary}\label{isolated}
Let $(M,\omega)$ be a closed symplectic manifold
equipped with a Hamiltonian circle action with a
moment map $\mu : M \rightarrow \R$. Assume that
fixed points are isolated. For each fixed point
$F \in M^{S^1}$, there exists a unique class
$\alpha_F \in H^{k_F}_{S^1}(M)$ such that
\begin{enumerate}
\item $\alpha_F |_{F'} = 0$ for every
$F' \in M^{S^1}$ with $\mu(F') < \mu(F)$,
\item $\alpha_F |_F = e^{-}_F = \prod w^{-}_iu$,
where $\{w^{-}_i \}_{1 \leq i \leq \frac{k_F}{2}}$
are the negative weights of
$S^1$-representation on $T_F M$, and
\item $\alpha_F |_{F'} = 0$
for every $F' \ne F \in M^{S^1}$ with $k_{F'}
\leq k_F$.
\end{enumerate}
Moreover, $\{ \alpha_F \}_{F \in M^{S^1}}$
is a basis of $H^*_{S^1}(M)$ as a $H^*(BS^1)$-module.
\end{corollary}

Now, consider a Hamiltonian $S^1$-manifold $(M,\omega)$
with a moment map $\mu : M \rightarrow \R$. Then there
exists an \textit{equivariant symplectic form}
$\widetilde{\omega}_\mu \in H^2_{S^1}(M)$ \textit{with
respect to $\mu,$} which satisfies the followings. (See \cite{Au} for the details.)

\begin{proposition}\label{easy}\cite{Au}
For a given closed Hamiltonian $S^1$-manifold with isolated fixed points with a moment map $\mu$, there exists an equivariant symplectic form $\widetilde{\omega}_\mu$ such that
\begin{enumerate}
\item $f^* \widetilde{\omega}_\mu = \omega,$ and
\item $[\widetilde{\omega}_\mu]|_F = -\mu(F)u$ for every $F \in M^{S^1}$.
\end{enumerate}

\end{proposition}

\section{Proof of Theorem \ref{theorem: main}} \label{section:
proof}

Let $(M,\omega)$ be a closed symplectic manifold with a
Hamiltonian circle action with only isolated fixed points.
Let $\mu : M \rightarrow \R$ be a moment map whose minimum
is zero and let $\Lambda_{2k} = \{ F_1^{2k}, \cdots,
F_{b_{2k}}^{2k} \}$ be the set of all fixed points of index
$2k$. Here, $b_{2k} := b_{2k}(M)$ is the $2k$-th Betti number of $M$.  Then $F_1^0$ is the unique minimum with $\mu(F_1^0) =
0$ by our assumption. Throughout this section, we assume
that $\mu$ is constant on each $\Lambda_{2k}$ so that
$\mu(F_1^{2k}) = \cdots = \mu(F_{b_{2k}}^{2k}) = c_{2k}$
for some $c_{2k} \in \R$ for every $k=0,1, \cdots,
n$. Let $\{\alpha_F ~| ~F \in M^{S^1} \}$ be the set of all canonical classes described in Corollary \ref{isolated}. To simplify our proof, we will denote by $\beta_F =
\frac{1}{\prod_i w^{-}_i} \alpha_F$, where $\{w^{-}_i \}_{1
\leq i \leq \frac{k_F}{2}}$ are the negative weights of
$S^1$-representation on $T_F M$ so that $\beta_F|_F =
u^{k_F}$ by Corollary \ref{isolated}.

\begin{lemma}\label{equiv symp}
For the equivariant symplectic class
$[\widetilde{\omega}_\mu] \in H^2_{S^1}(M)$, we have
$$[\widetilde{\omega}_\mu] = -c_{2k} \cdot (\beta_{F_1^2} +
\cdots + \beta_{F_{b_2}^2}).$$
\end{lemma}
\begin{proof}
Since $M$ is equivariantly formal by Theorem \ref{theorem:
formal}, the equivariant cohomology $H^*_{S^1}(M)$ is
isomorphic to a free $H^*(BS^1)$-module $H^*(M) \otimes
H^*(BS^1)$. Hence we have $H^2_{S^1}(M) \cong H^0(M)
\otimes H^2(BS^1) \oplus H^2(M) \otimes H^0(BS^1).$ By the
last statement of Corollary \ref{isolated}, $\{
u\beta_{F_1^0} = u, \beta_{F_1^2}, \cdots,
\beta_{F_{b_2}^2} \}$ forms a basis of $H^2_{S^1}(M)$. So
we may let $[\widetilde{\omega}] = a_0 \cdot u\beta_{F_1^0}
+ a_1 \cdot \beta_{F_1^2} + \cdots + a_{b_2} \cdot
\beta_{F_{b_2}^2}$ for some constants $a_0, a_1, \cdots,
a_{b_2}$.    Since each $\beta_{F_i^2}$ vanishes on $F_1^0$
for $i=1, \cdots, b_2$ by Corollary \ref{isolated}, we have
$$ [\widetilde{\omega}]|_{F_1^0} = a_0 \cdot u =
-\mu(F_1^0)u = 0 $$ by Proposition \ref{easy}. Hence we have $a_0 = 0$. Again, by
Corollary \ref{isolated}, each
$\beta_{F_i^2}$ vanishes on $F_j^2$ for every $j \neq i$ so
that we have
$$ [\widetilde{\omega}]|_{F_i^2} = a_i \cdot u =
-\mu(F_i^2)u = -c_{2k}u.$$ Therefore, we have $a_1 = \cdots
= a_{b_2} =c_{2k}.$ It finishes the proof.
\end{proof}

\begin{lemma}\label{vanish}
Let $\mu_{2k} = \mu - c_{2k}$ be the new moment map for
each $k=1, \cdots, n$. Then the equivariant symplectic
class $[\widetilde{\omega}_{\mu_{2k}}]$ vanishes on
$\Lambda_{2k}$.
\end{lemma}
\begin{proof}
By Proposition \ref{easy}, we have
$[\widetilde{\omega}_{\mu_{2k}}]|_{F} = -(\mu(F) - c_{2k})u^{\frac{k_F}{2}}$
for every fixed point $F \in M^{S^1}$. Hence we have
$$[\widetilde{\omega}_{\mu_{2k}}]_{F_i^{2k}} = -(\mu(F_i^{2k})
- c_{2k})u^k = 0$$ for every $i=1, \cdots, b_{2k}.$
\end{proof}

\begin{lemma}\label{distinct}
    The $n+1$ numbers $c_{2k}$'s are all distinct.
\end{lemma}

\begin{proof}
Assume that $c_{2i} = c_{2j}$ for some $i \neq j.$ Then the
class $$\eta =
[\widetilde{\omega}_{\mu_{0}}][\widetilde{\omega}_{\mu_{2}}]
\cdots \widehat{[\widetilde{\omega}_{\mu_{2i}}]}
\cdots[\widetilde{\omega}_{\mu_{2n}}] \in H^{2n}_{S^1}(M)$$
is an equivariant extension of $[\omega]^n \in H^{2n}(M)$. By Lemma \ref{vanish}, $\eta$ vanishes on
every fixed points $z \in \Lambda_{2k}$ for every $k$ so that $\eta$ is a zero class in $H^*_{S^1}(M)$. Since $[\omega]^n \neq 0$, it contradicts that $\eta$ is an extension of $[\omega]^n$.
\end{proof}

\begin{lemma}\label{zeroclass}
Let $\gamma \in H^{2k}_{S^1}(M)$ such that $\gamma|_z = 0$
for every $z \in \Lambda_{2i}$ with $i \leq k$. Then
$\gamma \equiv 0.$ Similarly, if $\gamma|_z = 0$ for every
$z \in \Lambda_{2i}$ with $i \geq n-k$, then $\gamma \equiv
0.$
\end{lemma}
\begin{proof}
By Theorem \ref{isolated}, the set $\{ u^{k-i} \cdot
\beta_{F_j^{2i}} \}_{i,j}$ is an $\R$-basis of
$H^{2k}_{S^1}(M)$ so that $\gamma = \sum_{i,j}
p_{i,j}u^{k-i} \cdot \beta_{F_j^{2i}}$ for some
coefficients $p_{i,j}$'s in $\R$. By our assumption, $\gamma$
vanishes on every fixed point of index less than or equal
to $2k$. Firstly, the restriction of $\gamma$ on the fixed
point $F_1^0$ with index zero is zero by our assumption so
that $$ \gamma|_{F_1^0} = p_{0,1}u^k = 0.$$ Hence we have
$p_{0,1} = 0.$ For index-two fixed points $\{F_1^2, \cdots,
F_{b_2(M)}^2 \}$, consider the restriction of $\gamma$ on
$F_j^2$ for any $j$. Then we have $$\gamma|_{F_j^2} =
p_{1,j}u^k = 0.$$ Hence we have $p_{1,j} = 0$ for every $j
\leq b_2(M)$. The rest part of the proof can be done in
exactly the same way.
\end{proof}

Now, we are ready to prove our main theorem.

\begin{proof}[Proof of Theorem \ref{theorem: main}]
For some $2k < n$, assume that
    \begin{displaymath}
        \begin{array}{cccc}
            [\omega]^{n-2k} : & H^{2k}(M) & \longrightarrow & H^{2n-2k}(M)\\
            & \alpha & \mapsto & \alpha \wedge [\omega]^{n-2k}\\
    \end{array}
\end{displaymath}
has a kernel $\gamma(\neq 0) \in H^{2k}(M),$ i.e. $\gamma
\cdot [\omega]^{n-2k} = 0.$ By Theorem \ref{theorem:
formal}, $M$ is equivariantly formal so that there is an
equivariant extension $\widetilde{\gamma} \in H^k_{S^1}(M)$
of $\gamma$ such that $f^*(\widetilde{\gamma}) = \gamma$
where $f : M \hookrightarrow M \times_{S^1} ES^1$ is an
inclusion as a fiber. Since the kernel of $f^*$ is an ideal
$u \cdot H^*_{S^1}(M)$ by Theorem \ref{theorem: formal}, we may choose $\widetilde{\gamma}$
such that $\widetilde{\gamma}|_z = 0$ for any fixed point
$z$ of index less than $2k.$ For each $k$, we denote by
$\mu_{2k}$ a moment map for the given circle action such
that $\mu_{2k} = 0$ on $\Lambda_{2k}$. Then the class
$$\delta = \widetilde{\gamma} \cdot
[\widetilde{\omega}_{\mu_{2k}}]
[\widetilde{\omega}_{\mu_{2k+2}}] \cdots
[\widetilde{\omega}_{\mu_{2n-2k-2}}] \in
H^{2n-2k}_{S^1}(M)$$ is an equivariant extension of $\gamma
\cdot [\omega]^{n-2k}$ satisfying $\delta|_z = 0$ for every
fixed point z of index less than $2n-2k$. Since $\gamma
\cdot [\omega]^{n-2k} = 0$ by our assumption, we have
$\delta \in \ker{f^*}$. Consequently, we have $\delta|_z =
0$ for every fixed point $z$ of index less than or equal to
$2n-2k.$ Hence $\delta \equiv 0$ by Lemma \ref{zeroclass}.
Note that if $\widetilde{\gamma}|_z \neq 0$ for some fixed
point $z$ of index greater than or equal to $2n-2k$, then
$$\delta|_z = \gamma|_z \cdot (c_{2k} - \mu(z)) \cdots (c_{2n-2k-2} - \mu(z)) \neq 0$$ by Lemma \ref{distinct}. Since $\delta
\equiv 0$, we have $\widetilde{\gamma}|_z = 0$ for every
fixed point $z$ of index greater than or equal to $2n-2k$.
By Lemma \ref{zeroclass} again, we have $\widetilde{\gamma}
\equiv 0$ so that $\gamma = 0$. It finishes the proof.
\end{proof}

\section{Examples} \label{section: example}
In the section, we give examples we said in
Section \ref{section: introduction}. First, we
introduce a well-known example of a
six-dimensional coadjoint orbit of $\SU(3).$

\begin{figure}[ht]
\begin{center}
\begin{pspicture}(-2,-2)(2,2.5) \footnotesize
\pspolygon[fillstyle=solid,fillcolor=lightgray,
linewidth=1pt](1,0)(0,1)(-1,1)(-1,0)(0,-1)(1,-1)(1,0)

\psline[linewidth=1.5pt](1,0)(0,1)(-1,1)(-1,0)(0,-1)(1,-1)(1,0)
\psline[linewidth=1.5pt](1,0)(-1,0)
\psline[linewidth=1.5pt](0,1)(0,-1)
\psline[linewidth=1.5pt](1,-1)(-1,1)

\psline[linewidth=0.5pt,arrowsize=5pt]{->}(0,-2)(0,2)
\psline[linewidth=0.5pt,arrowsize=5pt]{->}(-2,0)(2,0)
\psline[arrowsize=5pt]{->}(2,1)(1,2)

\psdots[dotsize=5pt](1,0)(0,1)(-1,1)(-1,0)(0,-1)(1,-1)

\uput[r](2,0){$x$} \uput[u](0,2){$y$}
\uput[ur](1,0){$1$} \uput[ur](0,1){$1$}
\uput[dl](-1,0){$-1$} \uput[dl](0,-1){$-1$}
\uput[dr](1,-1){$(1,-1)$}
\uput[ul](-1,1){$(-1,1)$}

\uput[ur](1.7,1.3){$\xi$}

\end{pspicture}
\end{center}
\caption{ \label{picture: hexagon} a moment map
image of a six-dimensional coadjoint orbit of
$SU(3)$}
\end{figure}
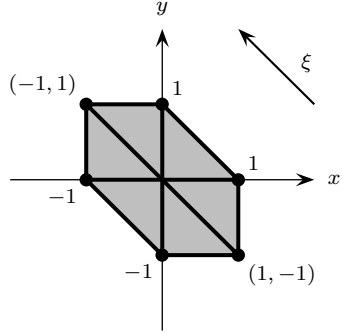

\begin{example}
The Lie algebra $\su(3)$ of $\SU(3)$ consists of
all traceless $3 \times 3$ skew-Hermitian
matrices. Let $T$ and $\torus$ be the standard
maximal torus of $\SU(3)$ and its Lie algebra,
i.e. $T$ and $\torus$ are subsets of diagonal
matrices in $\SU(3)$ and $\su(3),$ respectively.
We identify $\su(3)$ (resp. $\torus$) with its
dual $\su(3)^*$ (resp. $\torus^*$) through the
Killing form on $\su(3).$ Let $D_0 \in \torus^*$
be the diagonal element with entries $\sqrt{-1},$
$0,$ $-\sqrt{-1}.$ The \textit{coadjoint orbit}
$M \subset \su(3)^*$ of $\SU(3)$ through the
matrix $D_0$ is defined as the orbit of $D_0$ by
the coadjoint action of $\SU(3)$ on $\su(3)^*,$
i.e. conjugation. The coadjoint orbit endowed
with the Kostant-Kirillov symplectic form is a
symplectic manifold and (trivially) invariant
under the coadjoint $T$-action on $\su(3)^*.$ The
map
\[
\mu : M \rightarrow \torus^*, \quad (a_{ij})_{1
\le i,j \le 3} \longmapsto (\delta_{ij} \cdot
a_{ij})_{1 \le i,j \le 3}
\]
is a moment map for the $T$-action on $M,$ where
$\delta_{ij}$ is the Kronecker delta function.
Under the identification
\[
\torus^* \rightarrow \R^2, \quad
\left(
  \begin{array}{ccc}
    x \sqrt{-1}  & 0         & 0    \\
    0         & y \sqrt{-1}  & 0    \\
    0         & 0         & -(x+y) \sqrt{-1}  \\
  \end{array}
\right) \longmapsto (x,y),
\]
the moment map image is depicted in Figure
\ref{picture: hexagon}: black dots are images of
six fixed points. The vector $\xi=(-1 ,1)$ in the
Lie algebra $\mathfrak{t}$ defines a circle
subgroup in $T,$ and $\mu_\xi := \langle \mu, \xi
\rangle$ is a moment map for the circle action,
where $\langle~, \rangle$ is the evaluation
pairing between the dual space $\mathfrak{t}^*$
and $\mathfrak{t}.$ Then, we have
\begin{equation*}
\begin{array}{ll}
\mu_\xi = -2 ~ \text{ on } \Lambda_0, &
\quad \mu_\xi = -1 ~ \text{ on } \Lambda_2, \\
\mu_\xi = 1 ~ \text{ on } \Lambda_4, & \quad
\mu_\xi = 2 ~ \text{ on } \Lambda_6.
\end{array}
\end{equation*}
So, the moment map $\mu_\xi$ satisfies the
condition of the main theorem.
\end{example}

Second, we give a similar example of a
six-dimensional coadjoint orbit of $\SO(5).$

\begin{figure}[ht]
\begin{center}
\begin{pspicture}(-2,-2)(2,2.5) \footnotesize
\pspolygon[fillstyle=solid,fillcolor=lightgray,
linewidth=1pt](1,0)(0,1)(-1,0)(0,-1)(1,0)

\psline[linewidth=0.5pt,arrowsize=5pt]{->}(0,-2)(0,2)
\psline[linewidth=0.5pt,arrowsize=5pt]{->}(-2,0)(2,0)
\psline[arrowsize=5pt]{->}(4,-1)(3,2)
\psline[linewidth=1.5pt](1,0)(0,1)(-1,0)(0,-1)(1,0)
\psline[linewidth=1.5pt](-1,0)(1,0)
\psline[linewidth=1.5pt](0,-1)(0,1)

\psdots[dotsize=5pt](1,0)(0,1)(-1,0)(0,-1)

\uput[r](2,0){$x$} \uput[u](0,2){$y$}
\uput[ur](1,0){$1$} \uput[ur](0,1){$1$}
\uput[dl](-1,0){$-1$} \uput[dl](0,-1){$-1$}

\uput[ur](3.5,0.5){$\xi$}

\end{pspicture}
\end{center}
\caption{ \label{picture: square} a moment map
image of a six-dimensional coadjoint orbit of
$SO(5)$}
\end{figure}

\begin{example}
The Lie algebra $\so(5)$ of $\SO(5)$ consists of
$5 \times 5$ skew-symmetric matrices. For real
numbers $\theta$ and $x,$ denote by $R(\theta)$
and $S(x)$ the following two $2 \times 2$
matrices:
\[
\left(
  \begin{array}{cc}
    \cos \theta & -\sin \theta \\
    \sin \theta & \cos \theta \\
  \end{array}
\right) \quad \text{and} \quad \left(
  \begin{array}{cc}
    0 & x \\
    -x & 0 \\
  \end{array}
\right),
\]
respectively. Let $T$ and $\torus$ be the
standard maximal torus of $\SO(5)$ and its Lie
algebra, i.e.
\begin{align*}
T &= \Bigg\{ \left(
  \begin{array}{ccc}
    R(\theta)  & \mathbf{0}       & \mathbf{0}  \\
    \mathbf{0} & R(\theta^\prime) & \mathbf{0}  \\
    \mathbf{0} & \mathbf{0}       & 1     \\
  \end{array}
\right) \Bigg| ~ \theta, \theta^\prime \in \R
\Bigg\} \quad \text{and}
\\
\torus &= \Bigg\{ \left(
  \begin{array}{ccc}
    S(x)  & \mathbf{0}        & \mathbf{0}  \\
    \mathbf{0} & S(y)         & \mathbf{0}  \\
    \mathbf{0} & \mathbf{0}   & 0           \\
  \end{array}
\right) \Bigg| ~ x, y \in \R \Bigg\}.
\end{align*}
We identify $\so(5)$ and $\torus$ with their dual
$\so(5)^*$ and $\torus^*$ through the Killing
form on $\so(5),$ respectively. Let $D_0 \in
\torus^*$ be the element with $x=1$ and $y=0.$
For the coadjoint orbit $M \subset \so(5)^*$ of
$\SO(5)$ through the matrix $D_0$ and its
coadjoint $T$-action, the map
\[
\mu : M \rightarrow \torus^*, \quad (a_{ij})_{1
\le i,j \le 5} \longmapsto \left(
  \begin{array}{ccc}
    S(a_{12})  & \mathbf{0}   & \mathbf{0}  \\
    \mathbf{0} & S(a_{34})    & \mathbf{0}  \\
    \mathbf{0} & \mathbf{0}   & 0           \\
  \end{array}
\right)
\]
is a moment map for the $T$-action on $M.$ Under
the identification
\[
\torus^* \rightarrow \R^2, \quad \left(
  \begin{array}{ccc}
    S(x)  & \mathbf{0}        & \mathbf{0}  \\
    \mathbf{0} & S(y)         & \mathbf{0}  \\
    \mathbf{0} & \mathbf{0}   & 0           \\
  \end{array}
\right) \longmapsto (x,y),
\]
the moment map image is depicted in Figure
\ref{picture: square}. The vector $\xi=(-1,3)$ in
the Lie algebra $\mathfrak{t}$ defines a circle
subgroup in $T$ and $\mu_\xi := \langle \mu, \xi
\rangle$ is a moment map for the circle action.
Then, we have
\begin{equation*}
\begin{array}{ll}
(\mu_\xi +6)/2 = 0 ~ \text{ on } \Lambda_0, &
\quad (\mu_\xi +6)/2 = 2 ~ \text{ on } \Lambda_2, \\
(\mu_\xi +6)/2 = 4 ~ \text{ on } \Lambda_4, &
\quad (\mu_\xi +6)/2 = 6 ~ \text{ on } \Lambda_6.
\end{array}
\end{equation*}
So, the moment map $(\mu_\xi +6)/2$ is a
self-indexing morse function for the circle
action with respect to the Kostant-Kirillov
symplectic form multiplied by $1/2.$
\end{example}

Next examples are monotone semifree Hamiltonian
circle actions.

\begin{example}[Semifree Hamiltonian circle actions]\label{semifree}
Let $(M,\omega)$ be a $2n$-dimensional smooth compact
monotone symplectic manifold, i.e. a symplectic manifold
satisfying $\omega = c_1(M)$, where $c_1(M)$ is the first
Chern class of $M$ with respect to some $\omega$-tamed
almost complex structure $J$. Assume that $(M,\omega)$
admits a semifree Hamiltonian circle action with only
isolated fixed points. Then there exists a unique moment
map $\mu : M \rightarrow \R$ such that
$\widetilde{\omega}_{\mu} = c_1^{S^1}(M) \in
H^2_{S^1}(M;\R),$ where $c_1^{S^1}(M)$ is the equivariant
first Chern class of $M$. Note that the Morse index
$\lambda(z)$ is the twice number of negative weights of
tangential $S^1$-representation at $z$ for each fixed point
$z \in M^{S^1}$. Since we assumed that the action is
semifree, we have
$$-\mu(z) = i_z^*(\widetilde{\omega}_{\mu}) =
i_z^*(c_1^{S^1}(M)) = p_z - n_z = n - 2n_z$$ by
Proposition \ref{easy}, where $p_z$ ($n_z$,
respectively) is the number of positive (negative,
respectively) weights of tangential representation at
$z$. Hence we have $(\mu + n)(z) = 2n_z = \lambda(z)$
so that $\widetilde{\mu} := \mu + n$ is the
self-indexing moment map.
\end{example}

\begin{acknowledgement}
The second author is supported by GINUE research
fund.
\end{acknowledgement}

\end{document}